\documentclass[10pt]{article}

\usepackage{amsthm,amsmath,amssymb}

\usepackage{graphicx,url}

\usepackage[colorlinks=true,citecolor=blue,linkcolor=blue,urlcolor=blue]{hyperref}

\theoremstyle{plain}
\newtheorem{theorem}{Theorem}[section]
\newtheorem{lemma}[theorem]{Lemma}
\newtheorem{corollary}[theorem]{Corollary}

\theoremstyle{definition}
\newtheorem{definition}[theorem]{Definition}
\newtheorem{example}[theorem]{Example}

\theoremstyle{remark}

\title{\bf The $e-$positivity of some new classes of graphs }

\author{Stefan Mitrovi\'c\\
\small Faculty of Mathematics\\[-0.8ex]
\small University of Belgrade\\[-0.8ex]
\small Serbia\\
\small\tt stefan.mitrovic@matf.bg.ac.rs\\
\and
Tanja Stojadinovi\'c\\
\small Faculty of Mathematics\\[-0.8ex]
\small University of Belgrade\\[-0.8ex]
\small Serbia\\
\small\tt tanja.stojadinovic@matf.bg.ac.rs}

\begin{document}

\maketitle

\begin{abstract}
We introduce two classes of graphs - suns and dumbbells, both with few variations and explore their chromatic symmetric function and its $e$-positivity. We also give many connections of these two classes with other classes of connected graphs.

\bigskip\noindent \textbf{Keywords}: chromatic symmetric function, $e-$positivity, sun graphs, dumbbell graphs

\small \textbf{MSC2020}: 05E05, 05C15
\end{abstract}

\section{Introduction}

In 1995, Richard Stanley introduced the chromatic symmetric function as a generalization of the chromatic polynomial \cite{RS}. It inherits many properties of the chromatic polynomial, but not the property of deletion-contraction, a very useful tool for proving various identities for the chromatic polynomial using induction. Generalizations of the chromatic symmetric function, such as the chromatic quasisymmetric function \cite{SW} and the chromatic symmetric function in noncommuting variables \cite{GS}, allow for some form of a deletion-contraction relation. However, a triple-deletion property for ordinary chromatic symmetric function was discovered in \cite{OS}, and generalized to $k-$deletion in \cite{DvW}. 

Also, the chromatic symmetric function, being a more complicated combinatorial object, may distinguish some non-isomorphic graphs which cannot be distinguished by the chromatic polynomial. However, there are pairs of non-isomorphic graphs that have the same chromatic symmetric function, see \cite{MMW, OS, RS}. 

Nevertheless, the most active avenue of research is the positivity of the coefficients of the chromatic symmetric function of a given graph when expanded into elementary symmetric functions, known as $e-$positivity. A strong motivation for research in this direction is Stanley and Stembridge's conjecture posed in 1993, which states that the incomparability graph of $(3+1)-$ free poset is $e-$ positive \cite{SS}.

Wolfgang \cite{W} provided a very powerful criterion that any connected $e-$positive graph has a connected partition of every type, where a connected partition of a graph $G$ is a partition $\{V_1,\ldots, V_k\}$ of $V(G)$ such that each induced subgraph $G[V_i]$ is connected.

There are many graph classes that are known to be $e-$positive, such as complete graphs, paths, cycles, $K-$chains, lollipops, tadpoles and many more, \cite{DvW,GS,RS,WW}.

Further motivation for studying the positivity of chromatic symmetric functions is the close relationship between the Schur symmetric functions and representation theory. Namely, every $e-$positive graph is Schur positive since the coefficients of elementary symmetric functions in Schur basis are the Kostka numbers, which are nonnegative. Graphs that are shown to be Schur positive include the incomparability graphs of $(3+1)-$free posets, the incomparability graph of the natural unit interval order, edge 2-colorable hyperforests, some fan graphs, some broom and double broom graphs, and others \cite{G, G1, GP, SW, WW}.

 In this paper we introduce two classes of graphs, both with two variations. We examine their $e-$positivity by using Wolfgang's criterion and triple and $k-$deletion. Further, for some graphs that are not $e$-positive, we extract particular negative coefficient in the elementary basis. We also show how the problem of the $e-$ positivity of a large class of connected graphs can be reduced to that of one of these two classes. In addition, we investigate whether the chromatic symmetric polynomial or the chromatic symmetric function distinguish our graphs.

 Our paper is structured as follows. We cover the necessary background in Section 2. In Section 3, we discuss the $e-$ positivity of a new class of graphs that we call suns.  Further, in Subsection 3.1, we examine the positivity of some special types of suns in more detail. Section 4 introduces a class of $e-$positive graphs - dumbbells. We managed to deduce a recursive formula which allows us to inductively calculate their chromatic symmetric function in terms of the chromatic symmetric function of paths, cycles and complete graphs.

\section{Preliminaries}

We begin by recalling some algebraic, combinatorial and graph theoretic results that we will use throughout the paper. A \textit{partition} $\lambda=(\lambda_1, \lambda_2, \ldots, \lambda_{\ell(\lambda)})$ of $n$, denoted by $\lambda\vdash n$, is a list of positive integers whose parts $\lambda_i$ satisfy $\lambda_1\geq\lambda_2\geq\cdots\geq\lambda_{\ell(\lambda)}$ and $\sum_{i=1}^{\ell(\lambda)}\lambda_i=n$. The number $\ell(\lambda)$ is called the \textit{length} of $\lambda$. If $\lambda$ has exactly $m_i$ parts equal to $i$ for $1\leq i\leq n$, we may write $\lambda=(1^{m_1}, 2^{m_2}, \ldots, n^{m_n})$. By convention, there is one partition of $0$, denoted as $()$.

Related to every partition $\lambda=(\lambda_1, \lambda_2, \ldots, \lambda_{\ell(\lambda)})$ is its \textit{diagram}, which is a left-aligned table that, reading from top to bottom, has $\lambda_i$ cells in its $i$th row. By transposing this table, we get another partition's diagram and we call that partition the \textit{transpose} of $\lambda$, denoted as $\lambda^t$. Explicitly,
\[\lambda_i^t= \textrm{ number of parts of $\lambda\geq i$}.\]

Given two partitions, $\lambda, \mu\vdash n$, we will write $\lambda\leq_{p}\mu$ if the parts of $\mu$ are obtained by summing (not necessarily adjacent) parts of $\lambda$.

The \textit{algebra of symmetric functions} $\Lambda$ is a subalgebra of $\mathbb{Q}[[x_1, x_2, \ldots]]$ that consists of power sums of bounded degree that are invariant under the action of the symmetric group $S_n$, $n\in\mathbb{N}$, on the set of variables. This algebra has several natural bases.

The $i$th \textit{elementary symmetric function}, for $i\geq 1$, is given by 
\[e_i=\sum_{j_1<j_2<\cdots<j_i}x_{j_1}x_{j_2}\cdots x_{j_i},\]and for $i=0$ defined as $e_0:=1$. For a partition $\lambda=(\lambda_1, \lambda_2, \ldots, \lambda_{\ell(\lambda)})$, the \textit{elementary symmetric function} $e_{\lambda}$ is
\[e_{\lambda}=e_{\lambda_1}e_{\lambda_2}\cdots e_{\lambda_{\ell(\lambda)}}.\] The set $\{e_\lambda\}$ is a linear basis for $\Lambda$.

In many of our proofs, central place takes the \textit{power sum basis}. The $i$th \textit{power sum symmetric function} is defined as
\[p_0=1,\]
\[p_i=\sum_{j=1}^{\infty}x_j^i,\]
for $i\geq 1$. Analogously to the elementary case, we define 

\[p_{\lambda}=p_{\lambda_1}p_{\lambda_2}\cdots p_{\lambda_{\ell(\lambda)}}.\]

The $i$th power sum symmetric function is related to elementary functions in the following way \cite{DFvW}

\begin{equation}\label{piemu}
    p_i=\sum_{\mu=(1^{m1}, 2^{m_2}, \ldots, i^{m_i})}(-1)^{i-\ell(\mu)}\frac{i(\ell(\mu)-1)!}{\prod_{j=1}^{i}m_j!}e_{\mu}.
\end{equation}

    This relation means that when we calculate the coefficient of $e_{\mu}$ in a symmetric function written in the basis of power sum symmetric functions, we only need to focus on those $p_{\lambda}$ where $\mu\leq_{p}\lambda$.

Finally, our third basis of interest is the \textit{Schur basis}. Schur functions can be defined in different ways and have many beautiful combinatorial properties. For a partition $\lambda$, we will define them as

\[s_{\lambda}=\mathrm{det}(e_{\lambda_i^t -i+j})_{1\leq i, j\leq \lambda_1},\]
where we take $e_{\lambda_i^t -i+j}=0$ if $\lambda_i^t -i+j<0$.

If $\{u_\lambda\}$ is a basis of $\Lambda$, we will say that a symmetric function $f$ is $u$-\textit{positive} if all the coefficients $[u_\lambda] f$ in expression of $f$ in the basis $\{u_\lambda\}$ are nonnegative. Since the coefficients of elementary symmetric functions in Schur basis are well-known Kostka numbers, which are nonnegative, we have that every $e$-positive symmetric function is also $s$-positive.

From now on, we will focus on \textit{finite simple graphs}, their proper colourings and symmetric functions associated with them. Given a graph $G=(V, E)$, a \textit{proper colouring of $G$} is every function 
\[k: V\longrightarrow \mathbb{N}\]
such that, if $u, v\in V$ are adjacent in $G$, then $k(u)\neq k(v)$. The elements of $\mathbb{N}$ are called the \textit{colours}. If $e$ is an edge of $G$, we will write $G\setminus e$ for a graph obtained from $G$ by deleting edge $e$ and $G/e$ for a graph obtained from $G$ by contracting $e$ and identifying endpoints of $e$.

\begin{definition}\textit{The chromatic polynomial of a graph $G$, denoted as $\chi_G(x)$, for $x\in\mathbb{N}$, is a number of proper colourings of $G$ with $x$ colours.}
\end{definition}


One of the most useful properties, when dealing with the chromatic polynomial, is the following \textit{deletion-contraction property} that allows us to carry out proofs based on induction on the number of edges of $G$.

\begin{lemma}
\label{deletion}{(Deletion-contraction)} If $e$ is an edge of $G=(V, E)$, then $\chi_G=\chi_{G\setminus e}-\chi_{G/e}$.
\end{lemma}

Using this property, it can be easily deduced that $\chi_G(x)$ is indeed a polynomial in $x$. Moreover, we can calculate chromatic polynomial of any given graph. For example, for a tree $T$ with $d$ vertices, $\chi_{T}(x)=x(x-1)^{d-1}$. From this, we see that any two trees with the same number of vertices have the same chromatic polynomial. Hence, we would like to define a more sophisticated combinatorial object that will be a generalization of the chromatic polynomial.

\begin{definition} [\cite{RS}] \textit{For a graph $G$ with a vertex set $V=\{v_1, v_2, \ldots, v_d\}$, the chromatic symmetric function is defined to be
\[X_G=\sum_k x_{k(v_1)}x_{k(v_2)}\cdots x_{k(v_d)},\]
where the sum is over all proper colourings $k$ of $G$.}
\end{definition}

It is obvious that $X_G$ lies in $\Lambda$, so we can discuss about the positivity of $X_G$ in various bases. Also, it is easily seen that $X_G(1^n)=\chi_G(n)$ for every $n\in\mathbb{N}_{0}$, where $1^n$ means that $x_1= \cdots =x_n=1$ and $x_{n+1}=x_{n+2}= \cdots=0$. However, this object lacks the property of deletion-contraction.

There are many combinatorial expansions of the chromatic symmetric function. The most useful for us will be the one dealing with all subsets of the set of edges $E$ of $G$. Namely, if $G=(V, E)$ and $S\subseteq E$, we will write $\lambda(S)$ for a partition whose parts are equal to the numbers of vertices in the connected components of $(V, S)$.

\begin{lemma}
\label{expansion}
    \cite[Theorem 2.5]{RS} For a graph $G$ with a vertex set $V$ and an edge set $E$, we have \[X_G=\sum_{S\subseteq E}(-1)^{|S|}p_{\lambda(S)}.\]
\end{lemma}

Moreover, Orellana and Scott \cite{OS} used this expansion to prove a formula that we can take as a replacement for deletion-contraction property.

\begin{lemma}
\label{triple}
    \cite[Theorem 3.1]{OS} (Triple-deletion) If $e_1, e_2, e_3$ are edges of $G$ that form a triangle, then\[X_G=X_{G\setminus \{e_1\}}+X_{G\setminus \{e_2\}}-X_{G\setminus \{e_1, e_2\}}.\]
\end{lemma}

We say that a graph $G$ is $u$-positive, for $\{u_\lambda\}$ basis of $\Lambda$, if its chromatic symmetric function is $u$-positive. As we have mentioned before, a graph that is $e$-positive is necessarily $s$-positive. However, there are examples of graphs that are $s$-positive, but are not $e$-positive, see \cite{WW}. Wolfgang proposed a powerful criterion that allows us to test whether a graph is $e$-positive and it will be vital in some of our proofs.

\begin{theorem}\label{Wolfgang}\cite[Proposition 1.3.3]{W} If a connected $d$-vertex graph $G$ is $e$-positive, then $G$ has a connected partition of type $\lambda$ for every partition $\lambda\vdash d$.    
\end{theorem}
Hence, in order to prove that a graph $G=(V, E)$ is not $e$-positive, it will be sufficient to find a partition $\lambda\vdash |V|$ such that $G$ is missing a connected partition of type $\lambda$.

For some types of graphs, the coefficients of their chromatic symmetric function in the elementary basis can be calculated explicitly. Clearly, for complete graphs \[X_{K_n}=n!e_n.\] 
On the other hand, for paths and cycles we have the following two expansions.

\begin{theorem}
\label{putevi}
    \cite[Theorem 3.2]{Wo} Let $\lambda=(1^{a_1}, 2^{a_2},\ldots, d^{a_d})$ be a partition of $d$. The coefficient $c_{\lambda}$ of $e_{\lambda}$ in the expansion of $X_{P_d}$ is given by 
    \[c_{\lambda}=\binom{a_1+\cdots+a_d}{a_1, \cdots, a_d}\prod_{j=1}^d(j-1)^{a_j}+\]\[\sum_{i\geq 1}\left (\binom{(a_1+\cdots+a_d)-1}{a_1, \ldots, a_i-1,\ldots, a_d}\left (\prod_{j=1,\\ j\neq i}^{d}(j-1)^{a_j}\right )(i-1)^{a_i-1}\right ).\]
\end{theorem}

\begin{theorem}
\label{ciklusi}
    \cite[Theorem 3.5]{Wo} Let $\lambda=(1^{a_1}, 2^{a_2},\ldots, d^{a_d})$ be a partition of $d$. The coefficient $r_{\lambda}$ of $e_{\lambda}$ in the expansion of $X_{C_d}$ is given by 
    \[r_{\lambda}=\sum_{i\geq 1}\binom{(a_1+\cdots+a_d)-1}{a_1, \ldots, a_i-1, \ldots, a_d} i \prod_{j=1}^{d}(j-1)^{a_j}.\]
\end{theorem}

For some of our proofs, we will need an analogue of $X_G$, which is the \textit{chromatic symmetric function in noncommuting variables}, $Y_G$, introduced by Gebhard and Sagan in \cite{GS}. This point of view is convenient since it allows us to keep track of the colour which some proper colouring $k$ assigns to each vertex.
Just like the bases for $\Lambda$ are naturally indexed by partitions of numbers, the bases for algebra of noncommutative symmetric functions are naturally indexed by partitions of sets. Let $\Pi_d$ denote the \textit{lattice of set partitions} $\pi$ of $\{1, 2, \ldots, d\}$=[d], ordered by refinement. We will write $\pi=B_1/B_2/\ldots/B_k$ if $B_i$ are the \textit{blocks} of $\pi$. The
\textit{meet} (greatest lower bound) of the elements $\pi$ and $\sigma$ is denoted by $\pi\land\sigma$. We use 0 to denote the unique minimal element, and 1 for the unique maximal element. If $\pi$ is a partition of [d], we can naturally assign it a partition $\lambda(\pi)\vdash d$ whose parts are exactly the sizes of the blocks of $\pi$.

For a partition $\pi$ of [d], the \textit{noncommutative monomial symmetric function} $m_{\pi}$ is defined as
\[m_{\pi}=\sum_{i_1, i_2, \ldots, i_d} x_{i_1}x_{i_2}\cdots x_{i_d},\]
where the sum is over all sequences $i_1, i_2, \ldots, i_d$ of positive integers such that $i_j=i_k$ if and only if $j$ and $k$ are in the same block of $\pi$.

Still, we will mainly be interested in another basis, the basis of \textit{noncommutative elementary symmetric functions} 
\[e_{\pi}=\sum_{\sigma: \sigma\land\pi=0}m_{\sigma}=\sum_{i_1, i_2, \ldots, i_d} x_{i_1}x_{i_2}\cdots x_{i_d},\]
where the second sum runs over all sequences $i_1, i_2, \ldots, i_d$ of positive integers such that $i_j\neq i_k$ if $j$ and $k$ are in the same block of $\pi$.

\begin{definition}\cite{GS}
\textit{    Let $G$ be a graph with vertices labeled $v_1, v_2, \ldots, v_d$ in a fixed order. The chromatic symmetric function of $G$ in noncommuting variables is defined as
    \[Y_G=\sum_k x_{k(v_1)}x_{k(v_2)}\cdots x_{k(v_n)},\]
    where the sum is over all proper colourings $k$ of $G$, but the $x_i$ are noncommuting variables.}
\end{definition}

Unlike the regular chromatic symmetric function in commuting variables, this one allows some form of deletion-contraction, which is a great computational benefit. Also, many properties and expansions that hold for $X_G$ can be likewise translated to $Y_G$.

Let $\sigma, \pi$ be partitions of [d] and let $i\in[d]$. We write
\[\sigma\equiv_i \pi \textrm{ iff }\lambda(\sigma)=\lambda(\pi) \textrm{ and } |B_{\sigma, i}|=|B_{\pi, i}|, \]
where $B_{\sigma, i}$ and $B_{\pi, i}$ are respectively blocks of $\sigma$ and $\pi$ containing $i$. We naturally extend this definition to elementary functions 
\[e_{\sigma}\equiv_i e_{\pi}\textrm{ iff }\sigma\equiv_i\pi.\]
If we let $(\sigma)$ and $e_{(\sigma)}$ denote the equivalence classes of these relations, we can write
\[\sum_{\sigma\in\Pi_d} c_{\sigma}e_{\sigma}\equiv_i \sum_{(\sigma)\in\Pi_d/\equiv_i} c_{(\sigma)}e_{(\sigma)},\]
where $c_{(\sigma)}$ is the sum of all the coefficients $c_{\pi}$ from the left hand side such that $\pi\equiv_i\sigma$. We will say that a graph $G$ is $(e)-$positive if all the $c_{(\pi)}$ are nonnegative for some labeling of $G$ and suitably chosen congruence. 

\begin{theorem}\cite{GS}
\label{YtoX}
    If $Y_G$ is $(e)$-positive, then $X_G$ is $e$-positive.
\end{theorem}

Let $G$ be a graph with labeled vertices $V(G)=\{v_1, v_2, \ldots, v_d\}$ and edge set $E(G)$ and let $K_m$ denote the complete graph on $m$ vertices. We define $G+K_m$ to be the graph with the vertex set
\[V(G+K_{m})=V(G)\cup\{v_{d+1}, \ldots, v_{d+m-1}\}\]
and the edge set
\[E(G+K_m)=E(G)\cup\{v_iv_j\mid i, j\in [d, d+m-1]\}.\]

\begin{theorem}\cite[Theorem 7.6]{GS}
\label{lanac}
    IF $Y_G$ is $(e)$-positive, then $Y_{G+K_m}$ is $(e)$-positive.
\end{theorem}

\section{Suns}

We now come to our first object of study, the chromatic symmetric function of \textit{sun graphs}. Given the cycle $C_n$ ($n\geq 3$) with vertex set $\{v_1, v_2, \ldots, v_n\}$ and edge set $E=\{v_iv_{i+1}\mid i\in [n-1]\}\cup \{v_nv_1\}$ and a list of $n$ disjoint paths $P_{\lambda_1}, P_{\lambda_2}, \ldots, P_{\lambda_n}$, we obtain the \textit{(ordinary) sun graph} $S(n; \lambda_1, \lambda_2, \ldots, \lambda_n)$ by connecting $v_i$ with an edge to a leaf of $P_{\lambda_i}$, for every $i\in[n]$. If we replace the cycle $C_n$ in the previous definition with the complete graph $K_n$, we will get the \textit{complete sun graph} $\overline{S}(n; \lambda_1, \lambda_2, \ldots, \lambda_n)$.

In both of these two cases, the paths are called \textit{rays}, while $C_n$ in  $S(n; \lambda_1, \lambda_2, \ldots, \lambda_n)$, as well as $K_n$ in $\overline{S}(n; \lambda_1, \lambda_2, \ldots, \lambda_n)$ are called the \textit{body} of the sun. Vertex $v_i$ from the body is the \textit{beginning} of the ray $P_{\lambda_i}$ (although it does not belong to the ray - this definition is more convenient for us) and all the vertices of the sun with degree 1 are the \textit{ends} of rays.

\begin{figure}[h]
\label{primer sunca}
\caption{$S(5; 1, 2, 1, 2, 3)$ and $\overline{S}(4; 2, 1, 2, 1).$}
\centering
\includegraphics[width=1\textwidth]{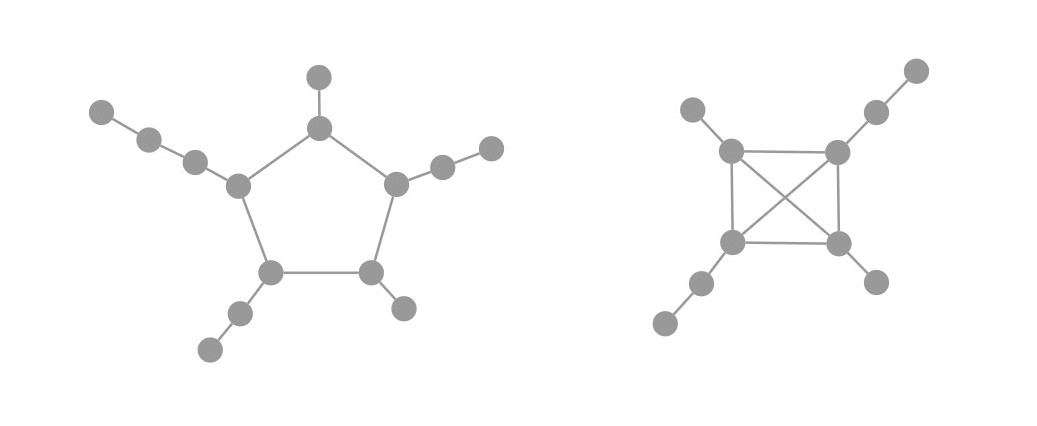}
\end{figure}

In \cite{FK}, our complete suns are called \textit{generalized spiders}, since they are exactly the \textit{line graphs} of ordinary spiders, where the line graph of $G=(V, E)$ is a new graph, denoted as $L_G$, whose vertex set is $E$ and two vertices $e_1, e_2$ in $L_G$ are adjacent if the edges $e_1, e_2$ share the same vertex in $G$. If $(\lambda_1, \lambda_2, \ldots, \lambda_d)$ is a partition, the spider $S(\lambda_1, \lambda_2, \ldots, \lambda_d)$ with $d$ legs is a graph consisting of $d$ paths $P_{\lambda_1}, P_{\lambda_2}, \ldots, P_{\lambda_d}$, each adjacent at a leaf to a center vertex. Many things are known about the $e$-positivity of spiders, see \cite{DSvW, WW, Z}. Another reason for investigating the chromatic symmetric function of suns is the possibility to reduce the $e$-positivity problem of a large class of connected graphs to examining $e$-positivity of suns. This is due to the fact that paths have all possible connected partitions, which gives us a solid ground for application of Wolfgang's criterion \ref{Wolfgang}.

\begin{lemma}
    Let $G_1, G_2, \ldots, G_n$ be connected graphs with $\lambda_1, \lambda_2, \ldots, \lambda_n$ vertices respectively and let $G$ be a graph obtained by connecting each vertex $v_i$ from the cycle $C_n=(\{v_1, v_2, \ldots, v_n\}, \{v_iv_{i+1}\mid i\in [n-1]\}\cup \{v_nv_1\})$ with an edge to a randomly chosen vertex $u_i$ from $G_i$. If $G$ has a connected partition of type $\mu$, then the sun $S(n; \lambda_1, \lambda_2, \ldots, \lambda_n)$ has a connected partition of type $\mu$ as well.
\end{lemma}

\begin{lemma}
\label{prebacivanje}
    Let $G_1, G_2, \ldots, G_n$ be connected graphs with $\lambda_1, \lambda_2, \ldots, \lambda_n$ vertices respectively and let $G$ be a graph obtained by connecting each vertex $v_i$ from the complete graph $K_n=(\{v_1, v_2, \ldots, v_n\}, \{v_iv_j\mid i\neq j\})$ with an edge to a randomly chosen vertex $u_i$ from $G_i$. If $G$ has a connected partition of type $\mu$, then the sun $\overline{S}(n; \lambda_1, \lambda_2, \ldots, \lambda_n)$ has a connected partition of type $\mu$ as well.
\end{lemma}

Let us dig deeper into finding connected partitions that certain suns are missing. For simplicity, when all the rays have the same size $k\geq 1$, we will write $S(n; k)$ and $\overline{S}(n; k)$ for the corresponding suns.

\begin{theorem}
\label{criterion}    The suns $S(n; k)$ and $\overline{S}(n; k)$ are missing a connected partition of certain type, thus they are not $e$-positive.
\end{theorem}

\begin{proof}
    Let us first consider the case $k=1$. Notice that graphs $S(n; 1)$ and $\overline{S}(n; 1)$ have $2n$ vertices. If $n$ is an even number, then both of these graphs are missing a connected partition of type $(n+1, n-1)$. This is because both parts of this partition are greater than $2$, so if we would assume that there is a connected partition of vertices of type $(n+1, n-1)$, the end and the beginning of any ray must be contained in the same part of this connected partition. This would imply that these parts are even, which is in contradiction with the fact that they are odd numbers. For the same reason, if $n$ is an odd number, we cannot construct a connected partition of type $(n, n)$. We can expand this idea to the case where $k>1$ as follows.

    If the number of vertices $n(k+1)$, where $k\geq 2$, is even, we can easily see that both $S(n; k)$ and $\overline{S}(n; k)$ do not have a connected partition of type $\lambda=(\frac{n(k+1)}{2}+1, \frac{n(k+1)}{2}-1)$. Suppose that such a connected partition exists. Then, since both parts of $\lambda$ are greater than $k$, the beginning and the end of any ray $P$ must be contained in the same part. Hence, the sizes of both parts must be divisible by $k+1$. This would imply that 2, which is exactly the difference of the parts of $\lambda$, is also divisible by $k+1\geq 3$, which is impossible.

    Similarly, if $n(k+1)$ is an odd number and $k\geq 2$, by using the same argument, one can deduce that $S(n; k)$ and $\overline{S}(n; k)$ are missing a connected partition of type $(\frac{n(k+1)+1}{2}, \frac{n(k+1)-1}{2})$, since the existence of such a connected partition would lead to the conclusion that $k+1\geq 3$ divides 1. 
\end{proof}

For the ordinary sun graphs, we can go a step further and calculate particular negative coefficients in the expansion of their chromatic symmetric function in the elementary basis. Although we know that $\overline{S}(n, k)$ are not $e$-positive as well, the same approach of finding a coefficient that is not positive can not be applied, since we would have deal with all of the connected spanning subgraphs of a complete graph.

\begin{lemma} \label{koeficijenti} For all $n\geq 3$, we have:

1) $[e_{(n+1, n-1)}]X_{S(n; 1)}=2n(1-n)$, if $n$ is even.

2) $[e_{(n, n)}]X_{S(n; 1)}=n(1-n)$, if $n$ is odd.

3) $[e_{(\frac{n(k+1)}{2}+1, \frac{n(k+1)}{2}-1)}] X_{S(n; k)}=n(k+1)(1-n)$, if $n(k+1)$ is even and $k\geq 2$.

4) $[e_{(\frac{n(k+1)+1}{2}, \frac{n(k+1)-1}{2})}]X_{S(n; k)}=n(k+1)(1-n)$, if $n(k+1)$ is odd and $k\geq 2$.  
\end{lemma}

\begin{proof}
    We will prove 1) and 4). Proofs of 2) and 3) are analogous. The idea is to calculate certain coefficients in the expansion of the chromatic symmetric function in the power sum basis, and later translate these results to results concerning the elementary basis.
    
    For 1), according to the comment immediately after Equation \ref{piemu}, we should only calculate $[p_{(n+1, n-1)}]X_{S(n; 1)}$ and $[p_{2n}]X_{S(n; 1)}.$ However, the expansion given in Lemma \ref{expansion} shows that the first of these two coefficients is 0 since $S(n; 1)$ is missing a connected partition of type $(n+1, n-1)$, as we have noted in the previous theorem. In order to calculate the second coefficient, we need to find all the subsets $S$ of the edge set $E$ of $S(n; 1)$ such that $S$ induces a connected spanning subgraph. Obviously, we can take $S=E$, where $|S|=2n$, which means that this subset contributes to the expansion from \ref{expansion} with the term $p_{2n}$. Further, if we take $|S|=|E|-1=2n-1$, we see that the removed edge has to be from the cycle since removing an edge from the rays would result in a disconnected graph. We can choose such an edge in $n$ ways, so these subsets give the term $-np_{2n}$. Removing two, or more edges from the sun results in a disconnected graph, hence, these results show that $[p_{2n}]X_{S(n; 1)}=1-n$. From the Equation \ref{piemu}, it follows that $[e_{(n+1, n-1)}]p_{2n}=2n$, and we finally get $[e_{(n+1, n-1)}]X_{S(n; 1)}=[e_{(n+1, n-1)}]p_{2n}\cdot [p_{2n}]X_{S(n; 1)}=2n(1-n)$.

    Similarly, for proving 4), we shall first focus on calculating coefficients $[p_{(\frac{n(k+1)+1}{2}, \frac{n(k+1)-1}{2})}]X_{S(n; k)}$ and $[p_{n(k+1)}]X_{S(n; k)}$. Just like in the previous part, the first coefficient is 0 and for the second one, we are searching for all the subsets $S$ of the edge set $E$ of $S(n; k)$ that induce a connected spanning subgraph. Since removing two or more edges gives a disconnected graph, there are only two possibilities for $|S|$. Firstly, for $|S|=|E|=n(k+1)$, we get the term $-p_{n(k+1)}$. Secondly, for $|S|=|E|-1$, we get $np_{n(k+1)}$, since the removed edge needs to be from the cycle. These two cases yield that $[p_{n(k+1)}]X_{S(n; k)}=n-1$. Finally, from the Equation \ref{piemu}, we conclude that $[e_{(\frac{n(k+1)+1}{2}, \frac{n(k+1)-1}{2})}]p_{n(k+1)}=-n(k+1)$, hence $[e_{(\frac{n(k+1)+1}{2}, \frac{n(k+1)-1}{2})}]X_{S(n; k)}=-n(k+1)(n-1)$.
\end{proof}

  We see that in all cases in the proof of Theorem \ref{criterion}, non-existence of a connected partition of certain type came from the fact that $n(k+1)$ can be written as a sum of two numbers that are not divisible by $k+1$.

\begin{theorem}
\label{gcd}
    If $\gcd(\lambda_1+1, \lambda_2+1, \ldots, \lambda_n+1)=m>1$, where $\lambda_1=\max \{\lambda_1, \lambda_2, \ldots, \lambda_n\}$ and $\lambda_1+2\leq\lambda_2+\cdots+\lambda_n+n-2$, then $S(n; \lambda_1, \lambda_2, \ldots, \lambda_n)$ and $\overline{S}(n; \lambda_1, \lambda_2, \ldots, \lambda_n)$ are missing a connected partition of type $\lambda=(\lambda_2+\cdots+\lambda_n+n-2, \lambda_1+2)$, hence are not $e$-positive.
\end{theorem}

\begin{proof}
    Assume that there is a connected partition of the vertex set of $S(n; \lambda_1, \lambda_2, \ldots, \lambda_n)$, or $\overline{S}(n; \lambda_1, \lambda_2, \ldots, \lambda_n)$, into sets $V_1$ and $V_2$ with $|V_1|=\lambda_2+\cdots+\lambda_n+n-2$ and $|V_2|= \lambda_1+2$. Since both parts of $\lambda$ are greater than $\lambda_i+1$ for every $i\in[n]$, we see that none of $V_1$ and $V_2$ can be contained in one ray only. This implies that if the beginning of some ray $P$ lies in $V_j$, for $j\in\{1, 2\}$, then $V_j$ certainly contains the whole ray $P$. We conclude that both parts of $\lambda$ are divisible by $m$. However, from $m\mid\lambda_1+1$, it follows that $m\nmid \lambda_1+2$ and similarly, from $m\mid \lambda_2+\cdots+\lambda_n+n-1$, we have that $m\nmid\lambda_2+\cdots+\lambda_n+n-2$.
\end{proof}

\begin{example}
    For $S(4; 5, 3, 3, 1)$ and $\overline{S}(4; 5, 3, 3, 1)$, $\lambda_i+1$ takes values 6, 4, 4, 2. Since $\gcd(6, 4, 4, 2)=2>1$ and $5+2<3+3+1+2$, previous theorem implies that these graphs are missing a connected partition of type $(9, 7)$.
\end{example}

We can also get a criterion for the $e$-positivity of these graphs by dealing with their \textit{perfect} and \textit{almost perfect matchings}, which are exactly the connected partitions of type $(2, 2, \ldots, 2)$ and $(2, 2, \ldots, 2, 1)$.  Since the end of every ray is of degree 1, there is only one way of pairing all the vertices from the rays. If a ray $P$ is of odd length, then its beginning must be paired with its neighbour from the ray, and if a ray $P$ is of even length, it must be paired with its neighbour from the body of the sun.

\begin{theorem}
  The sun $S(n; \lambda_1, \lambda_2, \ldots, \lambda_n)$ has a perfect (almost perfect) matching if and only if the induced subgraph on the set of vertices $\{v\mid v \textrm{ is the beginning of some ray of even length}\}$ has a perfect (almost perfect) matching.
\end{theorem}

Notice that the previous theorem does not include the complete suns. This is due to the fact that all the vertices from the body in a complete sun are adjacent, so the condition of the theorem always holds. But, when we consider ordinary suns, the fulfillment of this condition depends
on the arrangement of the rays.

\begin{example}
\label{razlika}
    $S(5; 3, 1, 2, 1, 2)$ does not have a perfect matching since two beginnings of rays of length 2 are not adjacent, but $S(5; 3, 2, 2, 1, 1)$ does.
\end{example}

We are able to apply the results we obtained for the suns to a large variety of graphs and improve Lemma \ref{prebacivanje}. In order to do so, we will need an obvious observation.

\begin{lemma}
    Let $G$ be a connected graph and let $G'$ be a graph obtained from $G$ by adding some edges. If $G'$ is missing a connected partition of type $\mu$, then $G$ is missing a connected partition of type $\mu$ as well.
\end{lemma}

Since any connected $n$-vertex graph can be obtained from the complete graph $K_n$ by removing some edges, we have the following lemma.

\begin{lemma}
\label{josprebacivanja}
     If $G_1, G_2, \ldots, G_n$ are connected graphs with $\lambda_1, \lambda_2, \ldots, \lambda_n$ vertices respectively and if $H$ is a connected graph with vertex set $V=\{v_1, \ldots, v_n\}$, let G be a graph obtained by connecting each vertex $v_i$ from $H$ with an edge to a randomly chosen vertex $u_i$ from $G_i$. If $G$ has a connected partition of type $\mu$, then the sun $\overline{S}(n; \lambda_1, \lambda_2, \ldots, \lambda_n)$ has a connected partition of type $\mu$ as well.
\end{lemma}

Hence, we could translate all the results we obtained so far that rely on finding a connected partition that some complete sun is missing, to many other graphs. For example, if instead of taking the cycle $C_n$, or the complete graph $K_n$ as the body of the sun, we take the spider $S(1^{n-1})$, we get the following corollary.

\begin{corollary}
    If the sun $\overline{S}(n; \lambda_1, \lambda_2, \ldots, \lambda_n)$ is missing a connected partition of type $\mu$, then, for any $i\in[n]$, a spider whose multiset of edges is $\{\lambda_j+1\mid j\neq i\}\cup\{\lambda_i\}$ is also missing a connected partition of type $\mu$ and is not $e$-positive.
\end{corollary}

\begin{example}
    According to Lemma \ref{prebacivanje}, the first one of the two graphs shown in Figure 2 is not $e$-positive since the sun $\overline{S}(5; 4)$ is missing a connected partition of some type. Lemma \ref{josprebacivanja} implies that the second one is not $e$-positive either.
\end{example}
    \begin{figure}[h]
    \label{pojednostavi}
\caption{}
\centering
\includegraphics[width=1\textwidth]{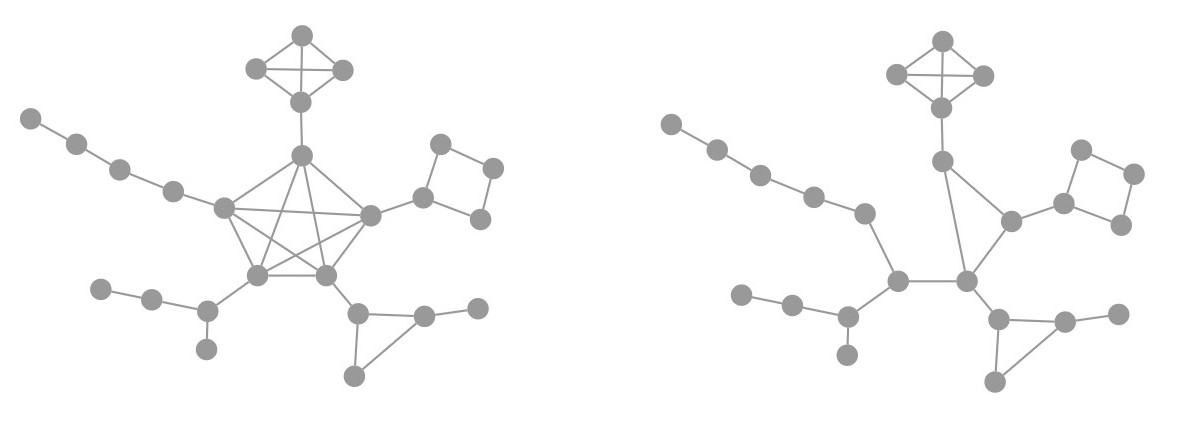}
\end{figure}

In \cite{FK}, it was proven that the chromatic symmetric functions of non-isomorphic complete suns are necessarily different. We could ask ourselves whether the same holds for ordinary suns. The thing that makes difference is that a multiset of integers $\{\lambda_1, \lambda_2, \ldots, \lambda_n\}$ uniquely determines the complete sun with rays of sizes $\lambda_i$, $i\in[n]$, while for the ordinary suns, we need to know additionally how the rays are arranged along the cycle, as we have seen in Example \ref{razlika}. Since $\chi_{C_n}(x)=(x-1)^{n}+(-1)^{n}(x-1)$ for any cycle $C_n$, $n\geq 3$, we easily see that $\chi_{S(n; \lambda_1, \lambda_2, \ldots, \lambda_n)}(x)=((x-1)^{n}+(-1)^{n}(x-1))(x-1)^{\lambda_1+\lambda_2+\cdots+\lambda_n}$. From this we have the following lemma.

\begin{lemma}
    If $S(n; \lambda_1, \lambda_2, \ldots, \lambda_n)$ and $S(m; \mu_1, \mu_2, \ldots, \mu_m)$ have the same chromatic symmetric functions, then $m=n$ and $\sum_{i=1}^n\lambda_i=\sum_{j=1}^n\mu_j$.
\end{lemma}

\subsection{The positivity of small suns}

In this section, we are going to specialize our discussion to the case of suns $S(3; \lambda_1, \lambda_2, \lambda_3)$. In this case, there is no difference between ordinary and complete suns. Particularly, when $\lambda_1=\lambda_2=\lambda_3=1$, we get \textit{the net} \cite{RS}, which is not $e$-positive. (Moreover, in \cite{FK} is proved that the \textit{generalized nets} are never $e$-positive.) 

On the other hand, in \cite{DSvW}, the authors proposed a conjecture that the only spiders that could be $e$-positive are exactly the spiders with three legs. Although this conjecture has not been proven yet, it is known that a spider with at least six legs can not be $e$-positive, see \cite{Z}. In \cite{Z}, another interesting conjecture is proposed - that a line graph of an $e$-positive spider is also $e$-positive. Notice that our suns $S(3, \lambda_1, \lambda_2, \lambda_3)$ are exactly the line graphs of spiders with three legs, such that all three legs have length greater than 1. The line graph of a three-legged spider, whose at least one leg is of length 1, is a chain of complete graphs, therefore $e$-positive according to Theorem \ref{lanac} . 

 Our first result concerning suns $S(3, \lambda_1, \lambda_2, \lambda_3)$, where the rays can have arbitrary length, is the simplification of the criterion given in Theorem \ref{gcd}.

\begin{theorem}
\label{pojednostav}
    If $\lambda_1=\max\{\lambda_1, \lambda_2, \lambda_3\}$ and $\lambda_1<\lambda_2+\lambda_3$, then the sun $S(3; \lambda_1, \lambda_2, \lambda_3)$ is missing a connected partition of type $(\lambda_2+\lambda_3+1, \lambda_1+2)$, thus is not $e$-positive.
\end{theorem}

\begin{proof}
   Suppose that there is a connected partition of the set of vertices of $S(3; \lambda_1, \lambda_2, \lambda_3)$ into two subsets $V_1$ and $V_2$, with $|V_1|=\lambda_2+\lambda_3+1$ and $|V_2|=\lambda_1+2$. Since both parts of this partition have sizes greater than $\lambda_i+1$ for $i\in \{1,2,3 \}$, none of them can be contained in one ray only. On the other hand, both parts of this partition have sizes smaller than $\lambda_i+\lambda_j+2$ for $i\neq j$, therefore, they can not contain more than one ray. 
\end{proof}

Similarly to the proof of Lemma \ref{koeficijenti}, it is easy to find some particular coefficients in the expansion of $X_{S(3; \lambda_1, \lambda_2, \lambda_3)}$ in the elementary basis.

\begin{lemma}
    If $\lambda_1=\max\{\lambda_1, \lambda_2, \lambda_3\}$ and $\lambda_1<\lambda_2+\lambda_3$, then:
    
    1) $[e_{(\lambda_2+\lambda_3+1, \lambda_1+2)}] X_{S(3; \lambda_1, \lambda_2, \lambda_3)}=-2(\lambda_1+\lambda_2+\lambda_3+3)$, if $\lambda_2+\lambda_3\neq\lambda_1+1$.

    2) $[e_{(\lambda_2+\lambda_3+1, \lambda_1+2)}] X_{S(3; \lambda_1, \lambda_2, \lambda_3)}=-(\lambda_1+\lambda_2+\lambda_3+3)$, if $\lambda_2+\lambda_3=\lambda_1+1$.
\end{lemma}
Notice that the conditions given in the previous lemma come from the fact that in case 2), the sizes of parts of the partition $(\lambda_2+\lambda_3+1, \lambda_1+2)$ are the same, while in case 1) they differ.

Consequently, from now on, we will only focus on suns $S(3; \lambda_1, \lambda_2, \lambda_3)$  with $\lambda_1\geq\lambda_2+\lambda_3$, where $\lambda_1=\max\{\lambda_1, \lambda_2, \lambda_3\}$. If we apply triple-deletion rule from Lemma \ref{triple} to $S(3; a, b, b)$, taking $e_1=uv$ and $e_2=uw$, with $u$ being the beginning of the ray of length $\lambda_1$ and $v, w$ being another two vertices from the body, we get
\[X_{S(3; a, b, b)}=2X_{S(a+1, b+1, b)}-X_{P_{2b+2}}X_{P_{a+1}},\]
where $S(m, n, p)$ represents the spider graph with three legs of length $m, n, p$. Since the paths are $e$-positive \cite{RS}, we have the following conclusion.

\begin{lemma}
    If $S(3; a, b, b)$ is $e$-positive, then $S(a+1, b+1, b)$ is also $e$-positive.
\end{lemma}

The $e$-positivity of spiders has recently attracted attention of many researchers. For example, in \cite{WW}, the authors proved that a spider $S(a, 2, 1)$ is $e$-positive if and only if $a\in\{3, 6\}$ and that a spider $S(a, 3, 2)$ is $e$-positive if and only if $a=5$. Therefore, according to the previous lemma, out of all suns $S(3; a, b, b)$ with $b\in\{1, 2\}$, only $S(3; 2, 1, 1)$, $S(3; 5, 1, 1)$ and $S(3; 4, 2, 2)$ are eligible for $e-$positivity. Nevertheless, an explicit computation of their chromatic symmetric functions shows that $[e_{(4,3)}]X_{S(3; 2, 1, 1)}=-2$ and $[e_{(4, 3, 3)}]X_{S(3; 5, 1, 1)}=-22$.

\begin{theorem}
\label{pocetak}
    The sun $S(a, b, b)$, for $b\in\{1, 2\}$, is $e$-positive if and only if $a=4$ and $b=2$.
\end{theorem}

In order to investigate the $e$-positivity of suns $S(3; a, b, b)$ for $b>2$, we will need a result from \cite{DSvW}.

\begin{lemma}
    \cite[Lemma 31]{DSvW} Suppose $S=S(\lambda_1, \lambda_2, \ldots, \lambda_d)$ is an $n$-vertex spider. Pick some $\lambda_i$ with $2\leq i<d$ and let $n=q(\lambda_i+1)+r$, where $0\leq r<\lambda_i+1$ and $t=\lambda_{i+1}+\cdots+\lambda_d$. If $q\geq \frac{\lambda_i+1}{t-1}$, then $S$ is not $e$-positive.
\end{lemma}

If we apply this result to the spider $S(a+1, b+1, b)$, taking $i=2$, we get the following:
\[q=\left \lfloor{\frac{a+2b+3}{b+2}}\right \rfloor = \left \lfloor{\frac{a-1}{b+2}+2}\right \rfloor=\left \lfloor{\frac{a-1}{b+2}}\right \rfloor+2.\]
Since $t=b$, previous lemma gives that if
\[\left \lfloor{\frac{a-1}{b+2}}\right \rfloor+2\geq \frac{b+2}{b-1},\]
then the spider $S(a+1, b+1, b)$ is not $e$-positive. This inequality can be rewritten as
\begin{equation}
    \left \lfloor{\frac{a-1}{b+2}}\right \rfloor\geq\frac{4-b}{b-1}.
\end{equation}
Since the left side of the last inequality is always nonnegative, we see that if $b\geq 4$, $S(a+1, b+1, b)$ is not $e$-positive. This implies that the same holds for the sun $S(3; a, b, b)$.

Hence, we only need to check what happens for $b=3$. We easily see that if $a\geq 6$, inequality from $(2)$ holds, therefore $S(a+1, b+1, b)$ and, consequently, $S(3; a, b, b)$, are not $e$-positive. Since we have already shown in Theorem \ref{pojednostav} that $a<6$ implies non $e$-positivity of $S(3; a, 3, 3)$, Theorem \ref{pocetak} gives the following beautiful result.

\begin{theorem}
    The sun $S(3; a, b, b)$ is $e$-positive if and only if $a=4$ and $b=2$.
\end{theorem}

There are several possible directions for further research of this subject. The natural question that arises at this point is whether there are some suns that are not $e-$positive, but are $s-$positive. Also, we could ask ourselves if there are infinitely many non-isomorphic $e-$positive suns.

\section{Dumbbells}

The second type of graphs we would like to introduce is closely related to \textit{tadpole} and \textit{lollipop} graphs. The tadpole graph $T_{m, l}$, for $m\geq 3,\ l\geq 1$, is obtained by connecting with an edge one vertex of the cycle $C_m$ and one leaf of the path $P_l$. If we take the complete graph $K_m$ instead of the cycle $C_m$ in the previous definition, we get the lollipop graph, denoted as $L_{m, l}$. We know that these graphs are $e$-positive and, in addition, we have the recursive relations for their chromatic symmetric functions \cite{DvW, WW1}.

The \textit{dumbbell} graph $D(m, l, n)$, for $m, n\geq 3$ and $l\geq 2$ is a graph that consists of two cycles $C_m$ and $C_n$ and the path $P_l$ such that one leaf of $P_l$ is connected by an edge to one vertex of $C_m$ and the other leaf to one vertex of $C_n$. In the previous definition, we can let $l=1$, and then our path is $P_1$, so the same vertex of $P_1$ is adjacent to both one vertex from $C_m$ and one vertex from $C_n$. The case where $l=0$ also makes sense, since we get two cycles connected by an edge. Moreover, we can introduce these graphs even for $l=-1$, as two cycles that share one vertex. If, instead of cycles $C_m$ and $C_n$, we take complete graphs $K_m$ and $K_n$ respectively, we get a graph that we will call the \textit{complete dumbbell} graph, $\overline{D}(m, l, n)$.

\begin{figure}[h]
\caption{$D(6, 1, 4)$ and $\overline{D}(5, 0, 4)$. }
\centering
\includegraphics[width=1\textwidth]{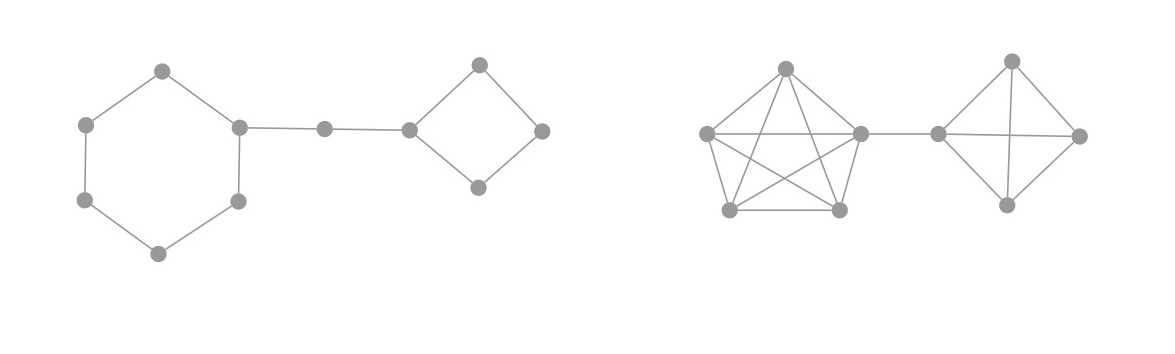}
\end{figure}

Notice that we have introduced a notation in which both $D(m, l, n)$ and $\overline{D}(m, l, n)$ have $m+l+n$ vertices. Obviously, $D(m, l, n)$ and $D(n, l, m)$, as well as $\overline{D}(m, l, n)$ and $\overline{D}(n, l, m)$ are pairwise isomorphic. In what follows, we are going to further explore the behaviour of the chromatic symmetric function of these graphs.

\subsection{The chromatic symmetric function of $D(m, l, n)$}

Recall that the chromatic polynomial of a cycle $C_n$ is $\chi_{C_n}(x)=(x-1)^n+(-1)^n(x-1)$. Likewise, the following formula holds
\[\chi_{D(m, l, n)}=\frac{1}{x}(x-1)^{l+3}((x-1)^{m-1}+(-1)^m)((x-1)^{n-1}+(-1)^n).\]

Therefore, if we know $\chi_{D(m, l, n)}$, we know the value of $l$ since $l+3$ is the greatest power of $(x-1)$ that divides the chromatic polynomial of a dumbbell. Further, if we apply the linear substitution $t=x-1$ in $P(x)=((x-1)^{m-1}+(-1)^m)((x-1)^{n-1}+(-1)^n)$, we get a polynomial in $t$, $Q(t)=t^{m+n-2}+(-1)^mt^{n-1}+(-1)^nt^{m-1}+(-1)^{m+n}$ and we are able to find the values of $m$ and $n$. Therefore, $\chi_{D(m, l, n)}$ contains enough information to uncover $m$, $l$ and $n$ from it. 

\begin{theorem}
    If $D(m, l, n)$ and $D(m', l', n')$ are non-isomorphic dumbbell graphs, then $X_{D(m, l, n)}\neq X_{D(m', l', n')}$.
\end{theorem}

In what follows, $P_2$ is identified  with $C_2$ and $C_n$ with $T_{n, 0}$ in favour of aesthetics of the formulas we obtained.

\begin{lemma}
\label{rekurzija}
    $X_{D(m, l, n)}=X_{D(m-1, l+1, n)}+X_{T_{n, m+l}}-X_{T_{n, l+1}}X_{C_{m-1}}$ for every $m>3$.
\end{lemma}

\begin{figure}[h]\caption{}\centering\includegraphics[width=1\textwidth]{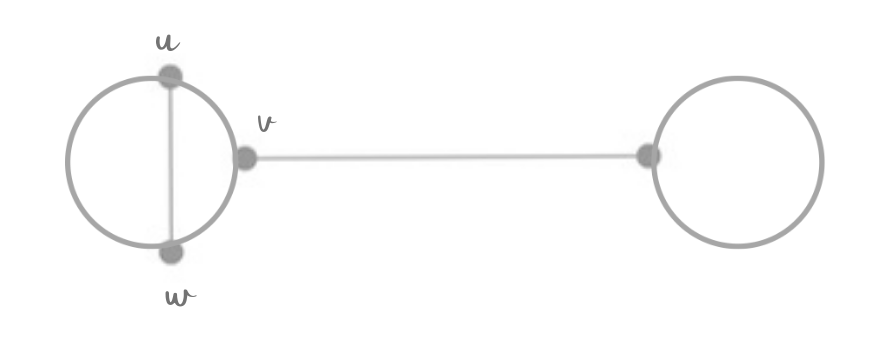}\end{figure}

\begin{proof}
    If $v$ is a vertex from the first cycle $C_m$ which has a degree greater than 2 and if $u$ and $w$ are the two neighbours of $v$ from the cycle, we will consider a graph $G$ obtained from ${D(m, l, n)}$ by adding an edge $uw$, see Figure 4. If $e_1=uv$, $e_2=vw$ and $e_3=wu$, triple-deletion from Lemma \ref{triple} gives

\[X_G=X_{G\setminus \{e_1\}}+X_{G\setminus\{e_2\}}-X_{G\setminus\{e_1, e_2\}},\]   \[X_G=X_{G\setminus \{e_2\}}+X_{G\setminus\{e_3\}}-X_{G\setminus\{e_2, e_3\}},\]

\noindent and, by combining these two we get \[X_{G\setminus\{e_3\}}=X_{G\setminus\{e_1\}}-X_{G\setminus\{e_1, e_2\}}+X_{G\setminus\{e_2, e_3\}},\] which is exactly what is stated.
\end{proof}
Notice that from the previous lemma we can not directly prove $e$-positivity of these graphs by induction. Moreover, in \cite{DSvW1} is proved that if $G$ and $H$ are two graphs, then $X_G-X_H$ is either $0$, or not $e$-positive, so we cannot anticipate how $X_{T_{n, m+l}}-X_{T_{n, l+1}}X_{C_{m-1}}$ will affects these coefficients. Still, if one of the cycles is a triangle, we have the desired $e$-positivity.

\begin{theorem}
    $X_{D(m, l, 3)}$ is $e$-positive.
\end{theorem}

\begin{proof}
We know that the cycles are $(e)$-positive, see \cite{GS}. Since our graphs $D(m, l, 3)$ can be obtained from the cycle $C_m$ as follows\[D(m, l, 3)=C_m+K_2+\ldots+K_2+K_3,\]
where $K_2$ appears $l+1$ times, Theorem \ref{lanac} implies that $Y_{{D(m, l, 3)}}$ is $(e)$-positive. According to Theorem \ref{YtoX}, $X_{D(m, l, 3)}$ is $e$-positive as well.
\end{proof}

The approach from the previous theorem can be applied even if we replace the triangle with any other complete graph $K_n$. Therefore, we may define the \textit{semicomplete dumbbell} graph $\Tilde{D}(m, l, n)$ as a graph that consists of the cycle $C_m$ and the complete graph $K_n$ that are linked by the path $P_l$. Analogously to the definition of the ordinary and the complete dumbbells, we may expand this definition to the case where $l=0$, or $l=-1$. Also, we can easily calculate that\[\chi_{\Tilde{D}(m, l, n)}(x)=((x-1)^{m}+(-1)^{m}(x-1))(x-1)^{l+2}(x-2)(x-3)\cdots (x-(n-1)).\] 

\begin{theorem}
    The semicomplete dumbbell graphs $\Tilde{D}(m, l, n)$ are $e$-positive. Non-isomorphic semicomplete dumbbell graphs have different chromatic symmetric functions.
\end{theorem}

By using an expanded version of techniques we described for the chromatic symmetric function in nocommuting variables, the authors of \cite{vWFAVW} succeeded in proving:

\begin{theorem} \cite[Proposition 6.7]{vWFAVW}
    $X_{D(m, l, n)}$ is $e$-positive.
 
\end{theorem}

Lemma \ref{rekurzija} allows us to express $X_{D(m, l, n)}$ in terms of chromatic symmetric functions of some simpler graphs. 

\begin{lemma}
\label{izraz}
   $ X_{D(m, l, n)}=(m-1)X_{T_{n, m+l}}-\sum_{k=1}^{m-2}X_{T_{n, l+k}}X_{C_{m-k}}$.
\end{lemma}

\begin{proof}
    We prove this by induction on $m\geq 3$. When $m=3$, let $v, u$ and $w$ be the vertices that form the triangle $C_3$, with $v$ being the one with the greatest degree, just like in Figure 4. If we apply triple-deletion \ref{triple}, taking $e_1=uv$, $e_2=vw$ and $e_3=wu$, we get
    \[X_{D(3, l, n)}=2X_{T_{n, l+3}}-X_{T_{n, l+1}}X_{C_2},\]
    as desired. Now, assume that the result holds for some $m-1\geq3$. Then, by Lemma \ref{rekurzija}, we have that
    \[X_{D(m, l, n)}=X_{D(m-1, l+1, n)}+X_{T_{n, m+l}}-X_{T_{n, l+1}}X_{C_{m-1}}=\]
    \[(m-2)X_{T_{n, m+l}}-\sum_{k=1}^{m-3}X_{T_{n, l+1+k}}X_{C_{m-1-k}}+X_{T_{n, m+l}}-X_{T_{n, l+1}}X_{C_{m-1}}=\]
    \[(m-1)X_{T_{n, m+l}}-\sum_{k=1}^{m-2}X_{T_{n, l+k}}X_{C_{m-k}}\]
    
\end{proof}

The expression in the previous lemma can further be simplified. Specifically, we would like to express the chromatic function of dumbbells in terms of the chromatic symmetric functions of cycles and paths, since we know the exact $e$-coefficients of these graphs, see Theorems \ref{putevi} and \ref{ciklusi}. Hence, the terms we would like to eliminate in expression in Lemma \ref{izraz} are the $X_{T_{a, b}}$. In \cite{WW1} is proved that
\[X_{T_{a, b}}=(a-1)X_{P_{a+b}}-\sum_{i=2}^{a-1}X_{P_{a+b-i}}X_{C_i}.\]
If we perform appropriate substitutions for every tadpole graph in \ref{izraz}, we get the following result.

\begin{theorem}
    \[X_{D(m, l, n)}=(m-1)(n-1)X_{P_{m+l+n}}-(m-1)\sum_{i=2}^{n-1}X_{P_{m+l+n-i}}X_{C_i}-\]\[(n-1)\sum_{j=2}^{m-1}X_{P_{m+l+n-j}}X_{C_j}+\sum_{i=2}^{n-1}\sum_{j=2}^{m-1}X_{P_{m+l+n-i-j}}X_{C_i}X_{C_j}.\]
\end{theorem}

\subsection{The chromatic symmetric function of $\overline{D}(m, l, n)$}

Complete dumbbells are chains of complete graphs, and consequently they are $e$-positive. They are also indifference graphs, just like their relatives, lollipops. \textit{Indifference} graph is a graph constructed by assigning a real number to each vertex and connecting two vertices by an edge when their numbers are within one unit from each other. If $m\geq n\geq 3$, then \[\chi_{\overline{D}(m, l, n)}(x)=x(x-1)^{l+3}(x-2)^2\cdots (x-(n-1))^2(x-n)\cdots (x-(m-1)).\]
Hence, $m$ is the greatest integer such that $(x-(m-1))\mid\chi_{\overline{D}(m, l, n)}$ and $n$ is the greatest integer such that $(x-(n-1))^2\mid\chi_{\overline{D}(m, l, n)}$.
\begin{theorem}
    If $\overline{D}(m, l, n)$ and $\overline{D}(m', l', n')$ are non-isomorphic complete dumbbell graphs, then $X_{\overline{D}(m, l, n)}\neq X_{\overline{D}(m', l', n')}$.
\end{theorem}
In the sequel, we will identify $\overline{D}(m, l, n)$ for $m=0, 1, 2$ with lollipops $L_{n, l+m}$.

\begin{lemma}
\label{rekurzija2}
    For every $m\geq 3$, we have that
    \[X_{\overline{D}(m, l, n)}=(m-1)X_{\overline{D}(m-1, l+1, n)}-(m-2)X_{K_{m-1}}X_{L_{n, l+1}}\]
\end{lemma}

\begin{figure}[h]\caption{}\centering\includegraphics[width=1\textwidth]{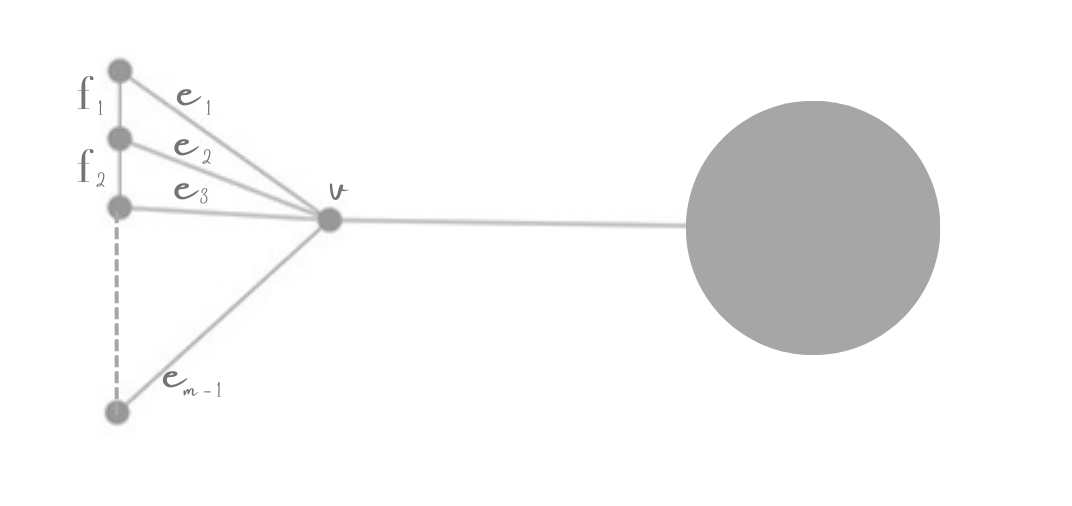}\end{figure}

\begin{proof}
    We will focus on the vertex $v$ of degree at least $m$ from the copy of $K_m$ in $\overline{D}(m, l, n)$. Out of all edges incident with $v$, there are $m-1$ edges that belong to the copy of $K_m$. Label these edges with $e_1, e_2, \ldots, e_{m-1}$. Also, let $f_i$ be the edge, which together with $e_i$ and $e_{i+1}$ forms a triangle. Note that if we remove from our graph any two subsets $S, T\subseteq \{e_1, e_2, \ldots, e_{m-1}\}$ of equal size $k$, then obtained graphs are isomorphic and thus, without loss of generality, we will focus on graphs of the form $\overline{D}(m, l, n)\setminus S_k$, where $S_k=\{e_1, e_2, \ldots, e_k\}$. If we focus on the triangle formed by edges $e_k, e_{k+1}$ and $ f_k$ and use triple-deletion from Lemma \ref{triple}, we get
    \[X_{\overline{D}(m, l, n)\setminus S_{k-1}}=X_{(\overline{D}(m, l, n)\setminus S_{k-1})\setminus\{e_k\}}+X_{(\overline{D}(m, l, n)\setminus S_{k-1})\setminus\{e_{k+1}\}}\]\[-X_{(\overline{D}(m, l, n)\setminus S_{k-1})\setminus\{e_k, e_{k+1}\}}= 2X_{\overline{D}(m, l, n)\setminus S_{k}}-X_{\overline{D}(m, l, n)\setminus S_{k+1}}.\] 
    If we apply this continually, and note that $\overline{D}(m, l, n)\setminus S_0=\overline{D}(m, l, n)$, $\overline{D}(m, l, n)\setminus S_{m-2}=\overline{D}(m-1, l+1, n)$ and $\overline{D}(m, l, n)\setminus S_{m-1}=K_{m-1}\cup L_{n, l+1}$, we get
    \[X_{\overline{D}(m, l, n)}=2X_{\overline{D}(m, l, n)\setminus S_{1}}-X_{\overline{D}(m, l, n)\setminus S_{2}}=\]
    \[2(2X_{\overline{D}(m, l, n)\setminus S_{2}}-X_{\overline{D}(m, l, n)\setminus S_{3}})-X_{\overline{D}(m, l, n)\setminus S_{2}}=\]
    \[3X_{\overline{D}(m, l, n)\setminus S_{2}}-2X_{\overline{D}(m, l, n)\setminus S_{3}}=\]
    \[3(2X_{\overline{D}(m, l, n)\setminus S_{3}}-X_{\overline{D}(m, l, n)\setminus S_{4}})-2X_{\overline{D}(m, l, n)\setminus S_{3}}=\]
    \[4X_{\overline{D}(m, l, n)\setminus S_{3}}-3X_{\overline{D}(m, l, n)\setminus S_{4}}=\]
    \[\vdots\]
    \[=(k+1)X_{\overline{D}(m, l, n)\setminus S_{k}}-kX_{\overline{D}(m, l, n)\setminus S_{k+1}},\]
    which holds for all $k\in [m-2]$. In particular, for $k=m-2$, we have:
    \[X_{\overline{D}(m, l, n)}=(m-1)X_{\overline{D}(m, l, n)\setminus S_{m-2}}-(m-2)X_{\overline{D}(m, l, n)\setminus S_{m-1}}=\]
    \[(m-1)X_{\overline{D}(m-1, l+1, n)}-(m-2)X_{K_{m-1}}X_{L_{n, l+1}}.\]
\end{proof}

This lemma also allows us to inductively calculate the chromatic symmetric function of $\overline{D}(m, l, n)$ in terms of the chromatic symmetric functions of complete graphs and lollipops.

\begin{lemma}
\label{prekololipopa}
For every $m\geq 3$,
   \[ X_{\overline{D}(m, l, n)}=(m-1)!X_{L_{n, m+l}}-\sum_{k=1}^{m-2}\frac{\prod_{i=1}^{k+1} (m-i)X_{K_{m-i}} X_{L_{n, l+i}}}{(m-k)}.\]
\end{lemma}

\begin{proof}
    We prove this by induction on $m\geq 3$. First, for $m=3$, let $v$ be the unique vertex from the copy of $K_3$ of degree at least 3 and let $u$ and $w$ be its two neighbours from the copy of $K_3$. The application of Lemma \ref{triple} with $e_1=uv$, $e_2=vw$ and $e_3=wu$, gives exactly what is desired. Now, assume that the result holds for some $m-1\geq3$. Then, by Lemma \ref{rekurzija2}, we have that
    \[X_{\overline{D}(m, l, n)}=(m-1)X_{\overline{D}(m-1, l+1, n)}-(m-2)X_{K_{m-1}}X_{L_{n, l+1}}=\]\[(m-1)\left ((m-2)!X_{L_{n, m+l}}-\sum_{k=1}^{m-3}\frac{\prod_{i=1}^{k+1} (m-1-i)X_{K_{m-1-i}} X_{L_{n, l+1+i}}}{(m-1-k)}\right )-(m-2)X_{K_{m-1}}X_{L_{n, l+1}}.\]
    After renaming $k:=k+1$ and $i:=i+1$ and multiplying the terms in the bracket with $m-1$, the previous expression becomes
    \[(m-1)!X_{L_{n, m+l}}-\sum_{k=2}^{m-2}\frac{\prod_{i=1}^{k+1}(m-i)X_{K_{m-i}}X_{L_{n, l+i}}}{(m-k)}-(m-2)X_{K_{m-1}}X_{L_{n, l+1}}.\]
    But the last term here can be seen as $\frac{(m-1)(m-2)}{(m-1)}X_{K_{m-1}}X_{L_{n, l+1}}$, which completes the proof.
    \end{proof}
Proposition $9$ from \cite{DvW} states that $X_{L_{a, b}}$ can be expressed in terms of the chromatic symmetric functions of paths and complete graphs in the following way
\[X_{L_{a, b}}=(a-1)!\left (X_{P_{a+b}}-\sum_{i=1}^{a-2}\frac{a-i-1}{(a-i)!}X_{K_{m-i}}X_{P_{n+i}}\right).\]
If we incorporate this in  \ref{prekololipopa}, we will get the following theorem.

\begin{theorem}
    For every $m, n\geq 3$ and for every $l\geq -1$,
    \[\frac{1}{(m-1)!(n-1)!}X_{\overline{D}(m, l, n)}=X_{P_{m+l+n}}-\sum_{i=1}^{m-2}\frac{(m-i-1)}{(m-i)!}X_{K_{m-i}}X_{P_{n+l+i}} \]
    \[-\sum_{j=1}^{n-2}\frac{(n-j-1)}{(n-j)!}X_{K_{n-j}X_{P_{m+l+j}}}+\sum_{i=1}^{m-2}\sum_{j=1}^{n-2}\frac{(m-i-1)(n-j-1)}{(m-i)!(n-j)!}X_{K_{m-i}}X_{K_{n-j}}X_{P_{l+i+j}}.\]
\end{theorem}




\bibliographystyle{model1a-num-names}
\bibliography{<your-bib-database>}







\end{document}